# A Cornucopia of Pythagorean triangles

*Konstantine Zelator*
*Department of Mathematics*
*301 Thackeray Hall*
*139 University Place*
*University of Pittsburgh*
*Pittsburgh, PA 15260*
*U.S.A*

*Also: Konstantine Zelator*
*P.O. Box 4280*
*Pittsburgh, PA 15203*
*U.S.A*

*e-mail addresses: 1) Konstantine_zelator@yahoo.com*
*2) kzet159@pitt.edu*

# 1. Introduction

As with many works in mathematics, the origin of this article lies in a discussion with a colleague about a classroom question. In Figure 1, two circles are illustrated, centered at $C_1$ and $C_2$. The two circles have only one point of intersection $I$; that point $I$ being the point of tangency between the two circles; they two circles being externally tangential. The two circles have three tangents or tangent lines in common. These are the two congruent tangents $\overline{T_1T_2}$ and $\overline{T_1'}\,\overline{T_2'}$, as well as their third common tangent which is the line perpendicular to $\overline{C_1C_2}$ at $I$; and which passes through the midpoint $M$ of the line segment $\overline{T_1T_2}$. The aforementioned subject of discussion, was simply the calculation of the length $T_1T_2$ of the tangent $\overline{T_1T_2}$, in terms of the two radii $C_1I = R_1$ and $C_2I = R_2$. This, as it turns out, is a simple matter, and the answer (as we will see below) is $2\sqrt{R_1R_2}$. This then led to a number theoretic exploration by this author, which resulted in this work. If we look at Figure 1, we can identify sixteen right triangles, which are listed in Section 2.

Now, consider the case when the two radii $R_1$ and $R_2$ are integers. In Section 6, we give the precise conditions the radii $R_1$ and $R_2$ must satisfy in order that all the sixteen right triangles ( listed in Section 2) are actually Pythagorean. In Section 3, we offer some immediate geometric observations from Figure 1, and in Section 4 we compute the side lengths of these Pythagorean triangles.
In Section 5, we state three results form number theory ( we offer a proof for Result 3) ; including the well known parametric formulas which describe the entire family of Pythagorean triples ( Result 1). In Section 7, we list the sidelengths of the 16 Pythagorean triangles and in Section 8 we present a numerical example. Finally, in Section 9, we explain why the diagonal lengths $d_1 = C_1T_2$ and $d_2 = C_2T_1$; are always irrational numbers , when the 16 right triangles are Pythagorean, this immediately follows form well known existing number theory results.

**Notation: 1)** If $X$ and $Y$ are two points on the plane, we will denote by $\overline{XY}$ the straight line segment joining the two points, and by $XY$ the length of the line segment $\overline{XY}$.

**2)** We will denote by $\overleftrightarrow{XY}$ the full straight line that goes through the two points $X$ and $Y$; and by $\overrightarrow{YX}$, the half-line or ray which emanates from $Y$ ( or whose vertex is the point $Y$ ),and which only contains those points on ( of the full line $\overleftrightarrow{XY}$ ) which lie on the side of (the point $Y$ ) which contains $X$.

**3)** If $X, Y, Z$ are three points on the plane we will denote by $< XYZ$ the angle formed by the rays $\overrightarrow{YX}$ and $\overrightarrow{YZ}$.

**4)** If $a$ and $b$ are natural numbers (or positive integers), $\gcd(a,b)$ will denote the greatest common divisor of $a$ and $b$.

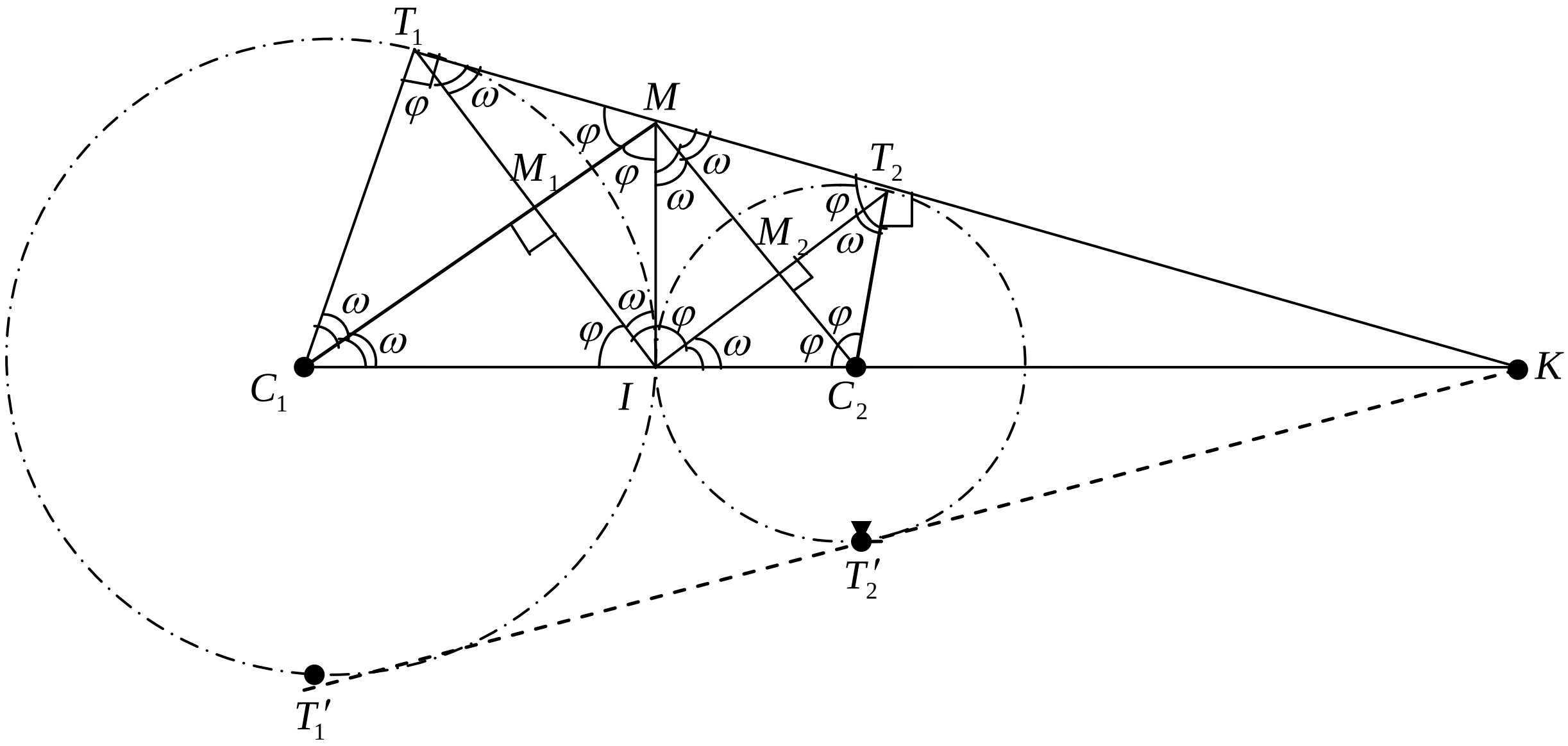

**Figure 1**

$\omega + \varphi = 90^\circ$, $C_1T_1 = C_1I = R_1$, $C_2I = C_2T_2 = R_2$ and $R_1 > R_2$

Also, from the right triangle $T_1IT_2$, since $M$ is the midpoint of $\overline{T_1T_2}$ ; we have

$$IM = T_1M = T_2M = \frac{T_1T_2}{2}$$

## 2. The sixteen right triangles

These are (from Figure 1) :

The two congruent right triangles $T_1M_1C_1$ and $C_1M_1I$ .

The two congruent right triangles $IM_2C_2$ and $C_2M_2T_2$.

The four congruent right triangles $T_1M_1M$ , $MM_1I$ , $IM_2M$ and $MM_2T_2$ .

The two congruent right triangles $C_1T_1M$ and $MIC_1$

The two congruent right triangles $C_2T_2M$ and $MIC_2$

The right triangle $C_1IM$

The right triangle $C_2IM$

The right triangle $C_1MC_2$

The right triangle $T_1IT_2$

The right triangle $C_2T_2K$

The right triangle $C_1T_1K$

## 3. Immediate geometric observations

Looking at Figure 1, we see that the geometry involved is pretty obvious and easy. First consider the various angels. We have,

$$< C_1T_1I =< T_1IC_1 =< T_1MC_1 =< C_1MI =< MIT_2 =< MT_2I =< IC_2M =< MC_2T_2 = \varphi$$
$$And\ < T_1C_1M =< IC_1M =< MT_1I =< T_1IM =< IMC_2 =< C_2MT_2 =< C_2IT_2 =< IT_2C_2 = \omega$$

And with $\omega + \varphi = 90^\circ$. Also, the lines $\overleftrightarrow{MI}$ and $\overleftrightarrow{C_1C_2}$ are perpendicular. Likewise, the lines $\overleftrightarrow{C_1M}$ and $\overleftrightarrow{MC_2}$ are perpendicular; and so are the lines $\overleftrightarrow{T_1I}$ and $\overleftrightarrow{IT_2}$ ; and the lines $\overleftrightarrow{MC_1}$ and $\overleftrightarrow{T_1I}$ are perpendicular as well. And both lines $\overleftrightarrow{T_1C_1}$ and $\overleftrightarrow{T_2C_2}$ are perpendicular to the line $\overleftrightarrow{T_1T_2}$ .

## 4. Length Computations

In this section we compute the side lengths of the 16 right triangles (listed in Section 2), in terms of the radii $R_1$ and $R_2$.
The key calculations are those of the lengths

$$T_1T_2\ , T_1I = a_1\ , T_2I = a_2\ , C_1M = x_1\ , and\ C_2M = x_2 .$$

The rest of the lengths follow easily from these five lengths.

Let $F$ be the foot of the perpendicular from the point $C_2$ to the line segment $\overline{C_1T_1}$ (Figure 2) Then ,

$$FC_1 = C_1T_1 - FT_1$$
$$= R_1 - C_2T_2$$
$$= R_1 - R_2$$

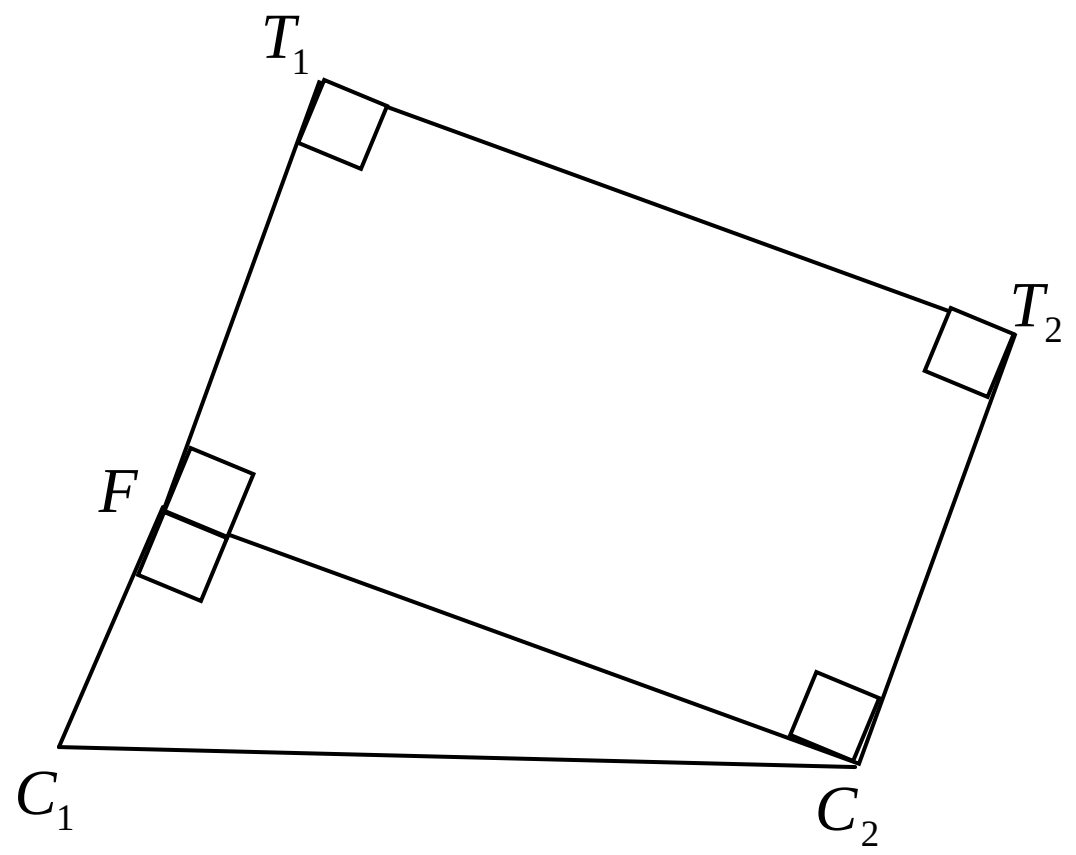

**Figure 2**

And from the right triangle $C_1FC_2$ we have,

$$(C_1C_2)^2 = (FC_1)^2 + (FC_2)^2$$
$$= (FC_1)^2 + (T_1T_2)^2;\ and\ by\ virtue$$

of $C_1C_2 = R_1 + R_2$ and $FC_1 = R_1 - R_2$ we obtain,

$$(T_1T_2)^2 = (R_1 + R_2)^2 - (R_1 - R_2)^2;$$
$$(T_1T_2)^2 = 4R_1R_2; \textit{or equivalently},$$

$$\boxed{T_1T_2 = 2\sqrt{R_1R_2}} \qquad (1)$$

From the right triangle $C_1IM$ we get

$$(C_1M)^2 = (C_1I)^2 + (IM)^2; \textit{and} \text{ since}$$
$$\text{IM} = \text{T}_1M = T_2M = \frac{T_1T_2}{2}, \textit{we obtain}$$
$$x_1^2 = R_1^2 + \left(\frac{T_1T_2}{2}\right)^2; \textit{and by} (1)$$
$$x_1^2 = R_1^2 + R_1R_2;$$

$$\boxed{x_1 = C_1M = \sqrt{R_1(R_1 + R_2)}} \qquad (2a)$$

Working similarly, from the right triangle $C_2IM$ we get

$$\boxed{x_2 = C_2M = \sqrt{R_2(R_1 + R_2)}} \qquad (2b)$$

Next, from the right triangles $C_1M_1I$ and $C_2M_2I$

we get $\left(\text{since } a_1 = T_1I, IM_1 = T_1M_1 = \frac{a_1}{2}, a_2 = T_2I, IM_2 = T_2M_2 = \frac{a_2}{2}\right)$,

$$\sin\omega = \frac{IM_1}{C_1I} = \frac{a_1/2}{R_1} = \frac{a_1}{2R_1} \text{ and } \sin\varphi = \frac{IM_2}{C_2I} = \frac{a_2}{2R_2}$$

So that, $\sin\omega = \dfrac{a_1}{2R_1}$ and $\sin\varphi = \dfrac{a_2}{2R_2}$ (3)

However, $\omega + \varphi = 90^\circ$ and so $\sin\omega = \cos\varphi$ and $\cos\omega = \sin\varphi$.

Therefore from the identity $\sin^2\omega + \cos^2\omega = 1$ and (3) we obtain

$$R_2^2 \cdot a_1^2 + R_1^2 \cdot a_2^2 = 4R_1^2R_2^2 \qquad (4)$$

Moreover from the right triangle $T_1IT_2$ we have,

$$(IT_1)^2 + (IT_2)^2 = (T_1T_2)^2; \text{ and by (1)}$$
$$a_1^2 + a_2^2 = 4R_1R_2 \qquad (5)$$

By solving the linear system of equations (4) and (5) in $a_1^2$ and $a_2^2$ (may use Kramer's Rule).

We further obtain

$$a_1^2 = \frac{4R_1^2 R_2 (R_1 - R_2)}{R_1^2 - R_2^2} = \frac{4R_1^2 R_2 (R_1 - R_2)}{(R_1 - R_2)(R_1 + R_2)} = \frac{4R_1^2 R_2}{R_1 + R_2}$$

Likewise, $a_2^2 = \dfrac{4R_2^2 R_1}{R_1 + R_2}$

Hence

$$\left\{ \begin{aligned} T_1 I &= a_1 = 2R_1 \cdot \sqrt{\frac{R_2}{R_1 + R_2}} \qquad (6a) \\ T_2 I &= a_2 = 2R_2 \cdot \sqrt{\frac{R_1}{R_1 + R_2}} \qquad (6b) \end{aligned} \right\}$$

And

$$\left\{ \begin{aligned} T_1 M_1 &= IM_1 = \frac{a_1}{2} = R_1 \cdot \sqrt{\frac{R_2}{R_1 + R_2}} \qquad (6c) \\ T_2 M_2 &= IM_2 = \frac{a_2}{2} = R_2 \cdot \sqrt{\frac{R_1}{R_1 + R_2}} \qquad (6d) \end{aligned} \right\}$$

Next we compute the height $h_1 = C_1 M_1$. From the right triangle $C_1 M_1 I$ we have,

$$(M_1 I)^2 + (C_1 M_1)^2 = (C_1 I)^2;$$

$$h_1^2 + \left(\frac{a_1}{2}\right)^2 = R_1^2; \; and\, by\, (6a)$$

$$h_1^2 + R_1^2 \cdot \left(\frac{R_2}{R_1 + R_2}\right) = R_1^2;$$

$$h_1^2 = R_1^2 \cdot \left[1 - \frac{R_2}{R_1 + R_2}\right] = \frac{R_1^3}{R_1 + R_2};$$

$$h_1 = R_1 \sqrt{\frac{R_1}{R_1 + R_2}}$$

Also $M_1 M = C_1 M - C_1 M_1 = x_1 - h_1 = \sqrt{R_1(R_1 + R_2)} - R_1 \sqrt{\dfrac{R_1}{R_1 + R_2}}$ .

Similarly we calculate $h_2 = C_2 M_2$ and $M_2 M = x_2 - h_2$ in terms of $R_1$ and $R_2$, from triangle $IM_2 C_2$.

Altogether we have

$$h_1 = C_1M_1 = R_1\sqrt{\frac{R_1}{R_1+R_2}} \qquad (7a)$$

$$M_1M = x_1 - h_1 = \sqrt{R_1(R_1+R_2)} - R_1\sqrt{\frac{R_1}{R_1+R_2}} \qquad (7b)$$

$$h_2 = C_2M_2 = R_2\sqrt{\frac{R_2}{R_1+R_2}} \qquad (7c)$$

$$M_2M = x_2 - h_2 = \sqrt{R_2\left(R_1+R_2\right)} - R_2\sqrt{\frac{R_2}{R_1+R_2}} \qquad (7d)$$

Furthermore, from the right triangle $T_1IT_2$ we have $IM = T_1M = T_2M = \frac{T_1T_2}{2}$ ; any by (1),

$$IM = T_1M = T_2M = \sqrt{R_1R_2} \qquad (8)$$

**Remark:** Note that $IM_1 = MM_2$ and $IM_2 = MM_1$, as it can be easily seen from (7b), (7d) , (6c), and (6d). This can also be seen from the fact that the quadrilateral $IM_1MM_2$ is a rectangle.

To finish the length computations; we must compute the lengths $C_2K$ , $C_1K$ , $T_2K$ and $T_1K$ .
From the similarity of the right triangles $C_2T_2K$ and $C_1T_1K$
we obtain

$$\frac{C_2K}{C_2T_2} = \frac{C_1K}{C_1T_1} \Leftrightarrow \frac{C_2K}{R_2} = \frac{C_2K+(R_1+R_2)}{R_1};$$ and solving for $C_2K$ yields

$$C_2K = \frac{R_2\left(R_1+R_2\right)}{R_1-R_2} \qquad (9a)$$

and thus, $C_1K = \frac{R_2\left(R_1+R_2\right)}{R_1-R_2} + R_1 + R_2$

or equivalently $C_1K = \frac{R_1\left(R_1+R_2\right)}{R_1-R_2}$ (9b)

We could compute $T_2K$ from the right triangles $C_2T_2K$ or alternatively, by using the similarity of the triangles $C_2T_2K$ and $C_1T_1K$ once more:

$$\frac{T_2K}{R_2} = \frac{T_2K+T_1T_2}{R_1}$$ ; and by (1)

we obtain, after solving for $T_2K$ ,

$$T_2K = \frac{2R_2\sqrt{R_1R_2}}{R_1 - R_2} \qquad (9c)$$

and thus, $T_1K = \dfrac{2R_2\sqrt{R_1R_2}}{R_1 - R_2} + 2\sqrt{R_1R_2}$

or equivalently, $T_1K = \dfrac{2R_1\sqrt{R_1R_2}}{R_1 - R_2}$ (9d)

## 5. Three results from number theory

The parametric formulas listed below in Results 1, describe the entire family of Pythagorean triples. A wealth of historical information on Pythagorean triangles can be found in references **[1]** and **[2]**.

**Result 1 :** *A triple* $(a,b,c)$ *of positive integers* $a,b,c$, *is a Pythagorean one with hypotenuse length* $c$ *(i.e.* $a^2 + b^2 = c^2$*) if , and only if (* $a$ *and* $b$ *may be switched),* $a = k(m^2 - n^2), b = k(2mn), c = k(m^2 + n^2)$, *for some positive integers* $k,m,n$ *such that* $m$ *and* $n$ *are relatively prime,* $m > n$, *and* $m,n$ *have different parities (i.e. one of them is even, the other odd) When* $k = 1$, *the Pythagorean triple is called primitive. Also, note that* $k = \gcd(a,b) = \gcd(a,c) = \gcd(b,c)$

The following result is well known and it can easily be found in number theory books; for example, see **[2]**.

**Result 2:** *Suppose that* $a,b$ *are relatively prime positive integers such that* $ab = c^2$, *where* $c$ *is a natural number. Then* $a = c_1^2$ *and* $b = c_2^2$, *where* $c_1, c_2$ *are relatively prime integers such that* $c_1c_2 = c$.

**Result 3:** *Suppose that* $a$ *and* $b$ *are positive integers such that* $ab = c^2$, *for some natural number* $c$. *Then,* $a = \delta r_1^2$ *and* $b = \delta r_2^2$, *for some positive integers* $\delta, r_1, r_2$*; such that* $r_1$ *and* $r_2$ *are relatively prime and* $\delta r_1 r_2 = c$.

**Proof :** Suppose that $ab = c^2$. *Let* $\delta$ *be the greatest common of* $a$ *and* $b$ *;* $\delta = \gcd(a,b)$ . Then $a = \delta c_1, b = \delta c_2$, for relatively prime positive integers $c_1$ and $c_2$. We obtain $\delta^2 c_1 c_2 = c^2$. Since $\delta^2$ is a divisor of $c^2$, it follows that $\delta$ is a divisor of $c$ ( more generally if $\delta^n$ is a divisor of $c^n$*;* then $\delta$ must be a divisor of $c$ - this is typically an exercise in a elementary number theory course). Put $c = \delta c_3$. We have, $\delta^2 c_1 c_2 = \delta^2 c_3^2; c_1 c_2 = c_3^2$ since $\gcd(c_1, c_2) = 1$, it follows from Result 2, that $c_1 = r_1^2$ and $c_2 = r_2^2$. And so, $r_1^2 r_2^2 = c_3^2$; $r_1 r_2 = c_3$.Altogether we have $a = \delta r_1^2$, $b = \delta r_2^2$, and $c = \delta r_1 r_2$, and with $r_1$ and $r_2$ being relatively prime □

Before we proceed further, note that by inspection we see from Figure 1 and because of the rectangle $IM_1MM_2$ that

$$x_1 = h_1 + \frac{a_2}{2} \; and \; x_2 = h_2 + \frac{a_1}{2}; where (from\, Section\; 4)\, x_1 = C_1M, h_1 = C_1M_1, x_2 = C_2M,$$

$$h_2 = C_2M_2, a_2 = IT_1, T_1M_1 = IM_1 = \frac{a_1}{2}, IM_2 = T_2M_2 = \frac{a_2}{2}.$$

## 6. The integrality of the lengths $T_1T_2, x_1, x_2$; and the rationality of the lengths $a_1, a_2, h_1, x_1 - h_1, h_2$ and $x_2 - h_2$

When $R_1$ and $R_2$ are integers, it is clear from (1) that $T_1T_2$ will be an integer precisely when the product $R_1R_2$ is a perfect or integer square. Otherwise the length $T_1T_2$ will be an irrational number. By Result 3, $R_1R_2$ will be an integer square if, and only if , $R_1 = \delta r_1^2$ and $R_2 = \delta r_2^2$, where $\delta, r_1, r_2$ are natural numbers with $r_1$ and $r_2$ being relatively prime. Now, if in addition to "the length $T_1T_2$ being an integer; we also require" the lengths $x_1$ and $x_2$ to be integers; it becomes apparent form formulas (2a) and (2b) , that since $R_1 = \delta r_1^2$ and $R_2 = \delta r_2^2$; $x_1$ and $x_2$ will be integers if, and only if, $r_1^2 + r_2^2 = r_3^2$, for some natural number $r_3$. In other words, precisely when $(r_1, r_2, r_3)$ is a primitive Pythagorean triple .

**Summary**

*When the two radii $R_1$ and $R_2$ are integers, then the three lengths $T_1T_2, x_1, x_2$, will be integers if, and only if, $R_1 = \delta r_1^2, R_2 = \delta r_2^2$, where $\delta, r_1, r_2$ are positive integers such that $r_1$ and $r_2$ are relatively prime and $r_1^2 + r_2^2 = r_3^2$, for some positive integer $r_3$. Then $(r_1, r_2, r_3)$ is a primitive Pythagorean triple and so, by Result 1,*

*Either $r_1 = m^2 - n^2$ and $r_2 = 2mn$; or alternatively $r_1 = 2mn$ and $r_2 = m^2 - n^2$*

*And in either case, $r_3 = m^2 + n^2$; where $m, n$ are relatively prime natural numbers such that $m > n$ and $m + n \equiv 1 (\mathrm{mod}\, 2)$ (i.e. one of $m, n$ is even, the other odd) And with $r_1 > r_2$, since $R_1 > R_2$*

Obviously, by (8) and (10), the lengths $IM = T_1M = MT_2$ are also integers. Furthermore, by (10), (6a)-(6d), (7a)-(7d), and (9a)-(9d), the ten lengths, $a_1, a_2, h_1, x_1 - h_1, h_2, x_2 - h_2, C_2K, C_1K, T_1K \text{ and } T_2K$ ; are rational numbers. Next, we calculate the above lengths in terms of the integers $\delta, r_1, r_2$ and $r_3$. The computations are straightforward. One simply uses (10) in conjunction with (6a)-(6d), (7a)-(7d), (8) and (9a)-(9d) in order to obtain the following:

$$T_1T_2 = 2\delta r_1 r_2 \quad (11a)$$

$$x_1 = \delta r_1 r_3 \quad (11b)$$

$$x_2 = \delta r_2 r_3 \quad (11c)$$

$$a_1 = \frac{2\delta r_2 r_1^2}{r_3} \quad (11d)$$

$$a_2 = \frac{2\delta r_1 r_2^2}{r_3} \quad (11e)$$

$$h_1 = \frac{\delta r_1^3}{r_3} \quad (11f)$$

$$h_2 = \frac{\delta r_2^3}{r_3} \quad (11g)$$

$$x_1 - h_1 = \delta r_1 r_3 - \frac{\delta r_1^3}{r_3} = \frac{\delta r_1 \cdot \left(r_3^2 - r_1^2\right)}{r_3} = \frac{\delta r_1 r_2^2}{r_3} = \frac{a_2}{2};$$
$$x_1 - h_1 = \frac{\delta r_1 r_2^2}{r_3} \quad (11h)$$

Similarly, $x_2 - h_2 = \dfrac{\delta r_2 r_1^2}{r_3}$ (11i)

$$IM = T_1M = T_2M = \delta r_1 r_2 \quad (11j)$$

$$C_2K = \frac{\delta r_2^2 r_3^2}{r_1^2 - r_2^2} \quad (11k)$$

$$T_2K = \frac{2\delta r_2^3 r_1}{r_1^2 - r_2^2} \quad (11l)$$

$$C_1K = \frac{\delta r_1^2 r_3^2}{r_1^2 - r_2^2} \quad (11m)$$

$$T_1K = \frac{2\delta r_1^3 r_2}{r_1^2 - r_2^2} \qquad \text{(11n)}$$

Now, since $(r_1, r_2, r_3)$ is a primitive Pythagorean triple; the three integers $r_1, r_2, r_3$ are pairwise relatively prime; $r_3$ is odd and one of $r_1, r_2$ is odd, the other even; and therefore $r_1^2 - r_2^2$ is also odd. The following coprimeness conditions follow readily (from (10) ) ( may also use $r_3^2 = r_1^2 + r_2^2$):

$$\left.\begin{array}{ll}
\text{In (11d)}, & \gcd(r_3, 2r_2r_1^2) = 1 \\
\text{In (11e)}, & \gcd(r_3, 2r_1r_2^2) = 1 \\
\text{In (11f)}, & \gcd(r_3, r_1^3) = 1 \\
\text{In(11g)}, & \gcd(r_3, r_2^3) = 1 \\
\text{In(11h)}, & \gcd(r_3, r_1r_2^2) = 1 \\
\text{In (11i)}, & \gcd(r_3, r_2r_1^2) = 1 \\
\text{In (11k)}, & \gcd(r_1^2 - r_2^2, r_2^2r_3^2) = 1 \\
\text{In (11l)}, & \gcd(r_1^2 - r_2^2, 2r_2^3r_1) = 1 \\
\text{In (11m)}, & \gcd(r_1^2 - r_2^2, r_1^2r_3^2) = 1 \\
\text{In (11n)}, & \gcd(r_1^2 - r_2^2, 2r_1^3r_2) = 1
\end{array}\right\} \qquad (12)$$

A fundamental result in number theory is Euclid's Lemma, which says that if an integer a divides the product $bc$ and $a$ is relatively prime to $b$, then $a$ must be a divisor of the integer $c$. Applying Euclid's Lemma to the formulas in (11d)-(11i) and (11k)-(11n) and in conjunction with the coprimeness conditions (12) the following becomes clear:

*In order that the lengths* $a_1, a_2, h_1, h_2, x_1 - h_1, x_2 - h_2, C_2K, KT_2, C_1K, KT_1$*; be integers as well, it is necessary and sufficient that the integer* $\delta$ *be divisible by both* $r_3$ *and* $r_1^2 - r_2^2$*; i.e.* $\delta$ *must be divisible by the least common multible of* $r_3$ *and* $r_1^2 - r_2^2$*. The least common multiple of* $r_3$ *and* $r_1^2 - r_2^2$ *is their product* $r_3(r_1^2 - r_2^2)$*, in virtue of the fact that* $\gcd(r_3, r_1^2 - r_2^2) = 1$.

Then we must have,

$$\delta = t \cdot r_3 \cdot \left(r_1^2 - r_2^2\right)$$
$$and\ (radii)\ R_1 = \delta r_1^2 = tr_3 r_1^2\left(r_1^2 - r_2^2\right)$$
$$and\ R_2 = \delta r_2^2 = tr_3 r_2^2\left(r_1^2 - r_2^2\right), \qquad (13)$$
$$where\, t\, is\, a\, positive\, \text{integer}$$

Consequently, from (11a)-(11n) and (13) we obtain the following length formulas:

$$T_1T_2 = 2tr_1r_2r_3\left(r_1^2 - r_2^2\right) \qquad (14a)$$

$$x_1 = tr_1r_3^2\left(r_1^2 - r_2^2\right) \qquad (14b)$$

$$x_2 = tr_2r_3^2\left(r_1^2 - r_2^2\right) \qquad (14c)$$

$$a_1 = 2tr_2r_1^2\left(r_1^2 - r_2^2\right) \qquad (14d)$$

$$a_2 = 2tr_1r_2^2\left(r_1^2 - r_2^2\right) \qquad (14e)$$

$$h_1 = tr_1^3\left(r_1^2 - r_2^2\right) \qquad (14f)$$

$$h_2 = tr_2^3\left(r_1^2 - r_2^2\right) \qquad (14g)$$

$$x_1 - h_1 = tr_1r_2^2\left(r_1^2 - r_2^2\right) \qquad (14h)$$

$$x_2 - h_2 = tr_2r_1^2\left(r_1^2 - r_2^2\right) \qquad (14i)$$

$$IM = T_1M = T_2M = tr_1r_2r_3\left(r_1^2 - r_2^2\right) \qquad (14j)$$

$$C_2K = tr_2^2r_3^3 \qquad (14k)$$

$$T_2K = 2tr_1r_3r_2^3 \qquad (14l)$$

$$C_1K = tr_1^2r_3^3 \qquad (14m)$$

$$T_1K = 2tr_2r_3r_1^3 \qquad (14n)$$

# 7. The sixteen Pythagorean triangles

Below we list the sixteen Pythagorean triangles obtained when all the lengths $T_1T_2, x_1, x_2, a_1, a_2, h_1, h_2, x_1-h_1, x_2-h_2, IM=T_1M=T_2M, C_2K, C_1K, T_1K, T_2K$ are integers. These lengths are calculated in terms of the integers $t, r_1, r_2, r_3$ as expressed in the formulas (14a)-(14n).

These 16 Pythagorean triangles are:

1) The two congruent Pythagorean triangles $C_1T_1M_1$ and $C_1M_1I$. They have hypotenuse length $C_1I=C_1T_1=R_1$, and leg lengths $C_1M_1=h_1$ and $T_1M_1=IM_1=\frac{a_1}{2}$.

2) The two congruent Pythagorean triangles $C_2T_2M_2$ and $C_2M_2I$. They have hypotenuse length $C_2I=C_2T_2=R_2$, and leg lengths $C_2M_2=h_2$ and $T_2M_2=IM_2=\frac{a_2}{2}$.

3) The four congruent Pythagorean triangles $T_1M_1M$, $MM_1I$, $IM_2M$, *and* $MM_2T_2$.

   The have hypotenuse length $T_1M=IM=T_2M=\frac{T_1T_2}{2}$; and leg lengths

   $T_1M_1=IM_1=MM_2=\frac{a_1}{2}=x_2-a_2$

   and also $MM_1=IM_2=T_2M_2=\frac{a_2}{2}=x_1-a_1$

4) The two congruent Pythagorean triangles $C_1T_1M$ and $C_1IM$. They have hypotenuse length $C_1M=a_1$ and leg lengths $C_1T_1=C_1I=R_1$ and $T_1M=IM=\frac{T_1T_2}{2}$.

5) The two congruent Pythagorean triangles $C_2T_2M$ and $C_2IM$. They have hypotenuse length $C_2M=a_2$ and leg lengths $C_2T_2=C_2I=R_2$ and $T_2M=IM=\frac{T_1T_2}{2}$.

6) The Pythagorean triangle $C_1MC_2$. It has hypotenuse length $C_1C_2=R_1+R_2$, and leg lengths $C_1M=x_1$ and $C_2M=x_2$.

7) The Pythagorean triangle $T_1IT_2$. It has hypotenuse length $T_1T_2$ and leg lengths $T_1I=a_1$ and $T_2I=a_2$.

8) The Pythagorean triangle $C_2T_2K$. It has hypotenuse length $C_2K$ and leg lengths $C_2T_2=R_2$ and $T_2K$.

9) The Pythagorean triangle $C_1T_1K$. It has hypotenuse length $C_1K$ and leg lengths $C_1T_1 = R_1$ and $T_1K$.

## 8. A numerical example

The first primitive Pythagorean triple $(r_1, r_2, r_3)$ with $r_1 > r_2$, is $(r_1, r_2, r_3) = (4,3,5)$. We take $t = 1$, and we apply formulas (13) and (14a)-(14n); in order to compute the numerical values of the various lengths. Specifically:

$$R_1 = 5 \cdot 4^2 \cdot (4^2 - 3^2) = 5 \cdot 16 \cdot 7 = 560$$
$$R_2 = 5 \cdot 3^2 \cdot (4^2 - 3^2) = 5 \cdot 9 \cdot 7 = 315$$
$$T_1T_2 = 2 \cdot 4 \cdot 3 \cdot 5 \cdot (4^2 - 3^2) = 2 \cdot 4 \cdot 3 \cdot 5 \cdot 7 = 840$$
$$x_1 = 4 \cdot 5^2 \cdot (4^2 - 3^2) = 4 \cdot 5^2 \cdot 7 = 700$$
$$x_2 = 3 \cdot 5^2 \cdot (4^2 - 3^2) = 3 \cdot 5^2 \cdot 7 = 525$$
$$a_1 = 2 \cdot 3 \cdot 4^2 \cdot (4^2 - 3^2) = 2 \cdot 3 \cdot 4^2 \cdot 7 = 672$$
$$a_2 = 2 \cdot 4 \cdot 3^2 \cdot (4^2 - 3^2) = 2 \cdot 4 \cdot 3^2 \cdot 7 = 504$$
$$h_1 = 4^3 \cdot (4^2 - 3^2) = 4^3 \cdot 7 = 448$$
$$h_2 = 3^3 \cdot (4^2 - 3^2) = 3^3 \cdot 7 = 189$$
$$x_1 - h_1 = 4 \cdot 3^2 \cdot (4^2 - 3^2) = 4 \cdot 3^2 \cdot 7 = 252$$
$$x_2 - h_2 = 3 \cdot 4^2 \cdot (4^2 - 3^2) = 3 \cdot 4^2 \cdot 7 = 336$$
$$IM = T_1M = T_2M = 4 \cdot 3 \cdot 5 \cdot (4^2 - 3^2) = 4 \cdot 3 \cdot 5 \cdot 7 = 420$$
$$C_2K = 3^2 \cdot 5^3 = 1125$$
$$T_2K = 2 \cdot 4 \cdot 5 \cdot 3^3 = 360$$
$$C_1K = 4^2 \cdot 5^3 = 2000$$
$$T_1K = 2 \cdot 3 \cdot 5 \cdot 4^3 = 1920$$

## 9. The irrationality of the diagonal lengths $d_1 = C_1T_2$ and $d_2 = C_2T_1$

There are two right triangles in Figure 1, that we have not mentioned thusfar. These are the triangles $C_1T_1T_2$ and $C_2T_2T_1$. As we have seen insofar, under the conditions (10) and (13); and the resulting formulas (14a)-(14n); sixteen Pythagorean triangles, listed in Section 7, are formed in Figure 1. The triangles $C_1T_1T_2$ and $C_2T_2T_1$ have integral leg lengths; these being $C_1T_1 = R_1$, $T_1T_2$, and $C_2T_2 = R_2$. However, as we will see below, the two hypotenuse lengths $d_1 = C_1T_2$ and $d_2 = C_2T_1$, are both irrational numbers (quadratic irrationals). Let us see why.

We have, $d_1^2 = R_1^2 + (T_1T_2)^2$; and by(10) and (11a) we obtain

$$d_1^2 = \delta^2 \cdot \left[r_1^4 + 4r_1^2 r_2^2\right];$$
$$d_1^2 = \delta^2 r_1^2 \cdot \left(r_1^2 + 4r_2^2\right); \qquad (15a)$$
$$d_1 = \delta r_1 \cdot \sqrt{r_1^2 + 4r_2^2}$$

And similarly, $d_2 = \delta r_2 \cdot \sqrt{4r_1^2 + r_2^2}$ (15b)

We point out that we do not need to make use of the special form that the integer $\delta$ has under (13). In other words, we will prove that both $d_1$ and $d_2$ are irrational under (10); and so when only the lengths $T_1T_2$, $x_1$, $x_2$ and $IM = T_1M = T_2M$ are integral for sure the rest of the lengths are rational, generally speaking, according to (11a)-(11n); the two lengths $d_1$ and $d_2$ must be irrational numbers. Indeed, consider (15a) and (10). We know that either $r_1 = m^2 - n^2$ and $r_2 = 2mn$ or vice-versa. Clearly $d_1$ will be either an integer or an irrational number, depending on whether $r_1^2 + 4r_2^2$ is an integer square or not. We have,

$$r_1^2 + 4r_2^2 = \left(m^2 - n^2\right)^2 + 4(2mn)^2;$$
$$r_1^2 + 4r_2^2 = m^4 + 14m^2n^2 + n^4$$

Recall from (10) that $m$ and $n$ are relatively prime integers, one being even, the other odd. Can it be,

$m^4 + 14m^2n^2 + n^4 = l^2$, for some positive integer $l$?

The answer is no. In reference **[1]**, Euler is mentioned as having proved that the diophantine equation $x^4 + 14x^2y^2 + y^4 = z^2$, with $\gcd(x, y) = 1$ can only have positive integer solutions when both $x$ and $y$ are odd.

In reference **[3]**, a proof can be found of the fact that the diophantine equation $x^4 + 14x^2y^2 + y^4 = z^2$ has no solutions in positive integers $x, y, z$ such that $\gcd(x, y) = 1$ and $x + y \equiv 1(\mathrm{mod}\, 2)$(i.e. one of $x, y$ is odd, the other even).Next, suppose that $r_1 = 2mn$ and $r_2 = m^2 - n^2$ We have $r_1^2 + 4r_2^2 = (2mn)^2 + 4\left(m^2 - n^2\right)^2 = 4\left(m^4 - m^2n^2 + n^4\right)$. Therefore , $r_1^2 + 4r_2^2$ will be a perfect square, precisely when $m^4 - m^2n^2 + n^4$ is a perfect square; which is impossible since the diophantine equation $x^4 - x^2y^2 + y^4 = z^2$ has a unique solution in positive integers $x, y, z$ such that $\gcd(x, y) = 1$. That solution is $x = y = z = 1$ (recall that $m$ and $n$ above have different parities). A proof of this result can be found in **[2]**. For the original paper that established the said result, refer to **[2]**. It is obvious from (15b), that the irrationality of $d_2$ is established in an identical manner.